\documentclass[a4paper]{amsart}
\usepackage{amssymb}
\usepackage{amscd}
\usepackage{amsthm}
\usepackage{amsmath}
\usepackage{mathtools}
\usepackage{latexsym}
\usepackage{hyperref}
\usepackage[all]{xy}
\usepackage[utf8]{inputenc}
\usepackage{enumitem}

\usepackage{tikz-cd}

\setlist[enumerate]{label={(\roman*)}}

\theoremstyle{plain}
\newtheorem{theorem}{Theorem}
\newtheorem{corollary}[theorem]{Corollary}
\newtheorem{lemma}[theorem]{Lemma}
\newtheorem{proposition}[theorem]{Proposition}

\theoremstyle{definition}
\newtheorem{definition}[theorem]{Definition}
\newtheorem{example}[theorem]{Example}

\theoremstyle{remark}

\newtheorem{remark}{Remark}

\numberwithin{theorem}{section}
\usepackage{newtxtext}
\usepackage{newtxmath}

\DeclareMathAlphabet\urwscr{U}{urwchancal}{m}{n}%
\DeclareMathAlphabet\rsfscr{U}{rsfso}{m}{n}
\DeclareMathAlphabet\euscr{U}{eus}{m}{n}
\DeclareFontEncoding{LS2}{}{}
\DeclareFontSubstitution{LS2}{stix}{m}{n}
\DeclareMathAlphabet\stixcal{LS2}{stixcal}{m} {n}
\newcommand{\Qcoh}[1]{\mbox{\rm Qcoh}(#1)}
\newcommand{\Spec}[1]{\operatorname{Spec}(#1)}

\newcommand{\Def}{\mbox{\rm{Def\,}}}

\newcommand{\im}{\mbox{\rm{Im\,}}}

\newcommand{\Hom}[3]{\operatorname{Hom}_{#1}(#2,#3)}
\newcommand{\Ext}[4]{\operatorname{Ext}^{#1}_{#2}(#3,#4)}
\newcommand{\Tor}[4]{\mbox{\rm{Tor}}_{#1}^{#2}(#3,#4)}

\newcommand{\rmod}[1]{\mbox{\rm{Mod}--}{#1}}

\newcommand{\lmod}[1]{{#1}\mbox{--\rm{Mod}}}
\newcommand{\ModR}{\text{Mod-}R}
\newcommand{\ModS}{\text{Mod-}S}

\begin{document}

\title{Flat relative Mittag-Leffler modules and Zariski locality}

\author{Asmae Ben Yassine and Jan Trlifaj}
\address{Charles University, Faculty of Mathematics
and Physics, Department of Algebra \\
Sokolovsk\'{a} 83, 186 75 Prague 8, Czech Republic}
\email{benyassine@karlin.mff.cuni.cz} 
\email{trlifaj@karlin.mff.cuni.cz}

\begin{abstract} The ascent and descent of the Mittag-Leffler property were instrumental in proving Zariski locality of the notion of an (infinite dimensional) vector bundle by Raynaud and Gruson in \cite{RG}. More recently, relative Mittag-Leffler modules were employed in the theory of (infinitely generated) tilting modules and the associated quasi-coherent sheaves, \cite{AH}, \cite{HST}. Here, we study the ascent and descent along flat and faithfully flat homomorphisms for relative versions of the Mittag-Leffler property. In particular, we prove the Zariski locality of the notion of a locally f-projective quasi-coherent sheaf for all schemes, and for each $n \geq 1$, of the notion of an $n$-Drinfeld vector bundle for all locally noetherian schemes.      
\end{abstract}

\date{\today}

\thanks{Research supported by GA\v CR 20-13778S. The first author also supported by Charles University Research Center program UNCE/SCI/022.}
	

\maketitle

\section{Introduction}

Relative Mittag-Leffler modules were introduced by Rothmaler in \cite{R1}. His approach was model theoretic: Mittag-Leffler modules were shown to be the counterparts of pure-injective modules in the sense that the former are atomic (i.e., they realize only the finitely generated pp-types) while the latter are saturated (i.e., they realize all pp-types). The adjective {\lq}relative{\rq} referred to restricting to theories of modules induced by definable subclasses of $\rmod R$. Much later, the important role of relative Mittag-Leffler modules for (infinite dimensional) tilting theory was recognized by Angeleri and Herbera \cite{AH}; this in turn led to a proof of finite type of all $1$-tilting modules in \cite{BH}. 

Flat Mittag-Leffler modules played a key role in proving Zariski locality of the notion of an (infinite dimensional) vector bundle in the classic work of Raynaud and Gruson, \cite[Seconde partie]{RG}. The locality follows by the Affine Communication Lemma (see e.g.\ \cite[5.3.2]{V}), whose assumptions are guaranteed by the ascent and descent of projectivity along flat ring homomorphisms, and faithfully flat ring homomorphisms, respectively.  

Once a structure theory of tilting modules over commutative rings was developed in \cite{APST} and \cite{HS}, it was possible to generalize the classic results to proving Zariski locality for various notions of quasi-coherent sheaves associated with tilting, \cite{HST}. Another generalization, employing the notion of a restricted flat Mittag-Leffler module, proved the Zariski locality of restricted Drinfeld vector bundles in \cite{EGT}. 

Our goal here is to refine the classic result on the ascent and descent of flat Mittag-Leffler modules to the relative setting. The main technical tools needed for this purpose are presented in Section \ref{tech}. In Section \ref{appl}, we apply these tools and prove Zariski locality of the corresponding notions of flat quasi-coherent sheaves. In particular, we prove the Zariski locality of the notion of a locally f-projective quasi-coherent sheaf for all schemes, and for each $n \geq 1$, of the notion of an $n$-Drinfeld vector bundle for all locally noetherian schemes. 

\section{Preliminaries}

Let $R$ be an (associative, unital) ring and $\rmod R$ the category of all (unitary right $R$-) modules. The elements of $\rmod R$ will often be referred to simply as \emph{modules}. Further, $\lmod R$ will denote the category of all (unitary) left $R$-modules. 

Let $n \geq 0$. A module $M$ is an \emph{FP$_n$ module} provided that $M$ possesses a projective resolution $\dots \to P_i \to P_{i-1} \to \dots \to P_0 \to M \to 0$ such that all the modules $P_i$ ($i \leq n$) are finitely generated. So FP$_0$ modules are just the finitely generated modules, FP$_1$ modules are the finitely presented ones, etc. Notice that the ring $R$ is right noetherian, iff the classes of FP$_n$ modules coincide for all $n \geq 0$ , while $R$ is right coherent, iff the classes of FP$_n$ modules coincide for all $n \geq 1$. 

We will denote by $\mathcal P _n$, $\mathcal F_n$, and $\mathcal I_n$ the classes of all modules of projective, weak, and injective dimension $\leq n$, respectively. 

Let $\mathcal B$ be a class of modules. Then ${}^\perp \mathcal B$ denotes the class of all modules $A$ such that $\Ext 1RAB = 0$ for each $B \in \mathcal B$. Similarly, $\mathcal B {}^\perp$ is the class of all modules $C$ such that $\Ext 1RBC = 0$ for all $B \in \mathcal B$. Further, $\mathcal B ^\intercal$ denotes the class of all left $R$-modules $D$ such that $\Tor 1RBD = 0$ for all $B \in \mathcal B$. Similarly, for a class of left $R$-modules $\mathcal D$, $^\intercal \mathcal D$ denotes the class of all modules $C$ such that $\Tor 1RCD = 0$ for all $D \in \mathcal D$. 

For a class of modules $\mathcal C$ we denote by $\varinjlim \mathcal C$ the class of all modules that are direct limits of direct systems consiting of modules from $\mathcal C$. For example, $\mathcal F_0 = \varinjlim \mathcal P_0$ for any ring $R$. Also, $\mathcal P \mathcal I$ will denote the class of all pure-injective modules.

We will need the following consequence of \cite[Theorem 8.40 and Corollary 8.42]{GT}:

\begin{lemma}\label{closure} Let $R$ be a ring and $\mathcal C$ be a class of FP$_2$-modules closed under extensions, direct summands and containing $R$. Let $\mathcal B = \mathcal C ^\perp$. Then $\varinjlim \mathcal C = {}^\perp (\mathcal B \cap \mathcal P \mathcal I)$. 
\end{lemma} 

We also recall the following identities satisfied by the Tor bifunctor.

\begin{lemma}\label{Ext-Tor} Let $\varphi : R \to S$ be a flat ring homomorphism of commutative rings.
\begin{enumerate}
\item[(1)] For all modules $A$ and $B$, there is an $S$-isomorphism $\Tor 1RAB \otimes _R S \cong \Tor 1S{A\otimes_R S}{B \otimes_R S}$.
\item[(2)] If $A$ is a module and $B$ is an $S$-module, then there is an $S$-isomorphism $\Tor 1S{A\otimes_R S}B \cong \Tor 1RAB$.
\end{enumerate}
\end{lemma}  
\begin{proof} (1) is a particular instance of \cite[Theorem 2.1.11]{EJ}, and (2) a particular instance of \cite[Proposition VI.4.1.2]{CE}.
\end{proof}

The central notion of our paper is that of a relative Mittag-Leffler module:

\medskip
\begin{definition}\label{mldef}
Let $R$ be an arbitrary ring, $M \in \rmod R$ and $\mathcal Q \subseteq \lmod R$. Then $M$ is \emph{$\mathcal Q$-Mittag-Leffler} (or \emph{Mittag-Leffler relative to $\mathcal Q$}), provided that the canonical morphism $\psi_M : M \otimes_{R} {\prod_{i \in I} Q_i} \to \prod_{i \in I} M \otimes_R Q_i$ defined by $\psi_M (m \otimes (q_i)_{i \in I}) = (m \otimes q_i )_{i \in I}$ is injective for any family $( Q_i \mid i \in I )$ consisting of elements of $\mathcal Q$. 
\end{definition}

As mentioned above, relative Mittag-Leffler modules were introduced in \cite{R1}. Further results on these modules were proved in \cite{AH}, and in the more recent papers \cite{H} and \cite{R2}. Following \cite{HT}, we will denote by $\mathcal D _{\mathcal Q}$ the class of all flat $\mathcal Q$-Mittag-Leffler modules. 

\medskip
The two borderline cases of Definition \ref{mldef} occur for $\mathcal Q = \emptyset$, when $\mathcal D _{\mathcal Q} = \mathcal F_0$, and for $\mathcal Q _0 = \lmod R$, when $\mathcal D _{\mathcal Q _0} = \mathcal F \mathcal M$ is the class of all flat Mittag-Leffler modules. As $\mathcal Q _0 = (\mathcal F _0)^\intercal$, the latter setting can be extended as follows: for each $n \geq 0$, we let $\mathcal Q _n = (\mathcal F _n)^\intercal$. We will call the modules $M \in \mathcal D _{\mathcal Q _n}$ \emph{flat $n$-Mittag-Leffler}. 

Another case of interest is when $\mathcal Q = \{ R \}$, or equivalently, $\mathcal Q$ is the class of all flat left $R$-modules. Then the flat $\mathcal Q$-Mittag-Leffler modules coincide with the \emph{f-projective} modules, that is, the modules $M$ such that each homomorphism from a finitely generated module to $M$ factorizes through a free module, see \cite{G} or \cite[\S 3]{BT}. 

Denoting the class of all f-projective modules by $\mathcal F \mathcal P$, we have the following chain of classes of modules 
$$(\ast) \quad \mathcal P _0 \subseteq \mathcal F \mathcal M = \mathcal D _{\mathcal Q _0} \subseteq \dots \subseteq \mathcal D _{\mathcal Q _{n}} \subseteq \mathcal D _{\mathcal Q _{n + 1}} \subseteq \dots \subseteq \mathcal F \mathcal P \subseteq \mathcal F _0.$$
\medskip

The inclusions in the chain $(\ast)$ need not be strict in general. For example, if $R$ has weak global dimension $\leq n$, then $\mathcal D _{\mathcal Q _{n}} = \mathcal D _{\mathcal Q _{n + 1}} = \dots = \mathcal F \mathcal P$. If $R$ is a right perfect ring, then all the classes in the chain $(\ast)$ coincide.   

\begin{remark}\label{remsim} 1. Other variants of the notion of a flat Mittag-Leffler module, called \emph{restricted flat Mittag-Leffler} modules, were introduced in \cite{EGT}. Their classes form a chain located between the classes $\mathcal P_0$ and $\mathcal F \mathcal M$. 

2. The following generalization of the notion of an f-projective module goes back to Simson \cite{S}: given a cardinal $\kappa \geq \aleph_0$, a module $M$ is \emph{$\kappa$-projective} if each homomorphism from a $< \kappa$-generated module to $M$ factorizes through a free module. Denote by $\mathcal C _\kappa$ the class of all $\kappa$-projective modules. Since $\mathcal C _{\aleph_0} = \mathcal F \mathcal P = \mathcal D _{\{ R \}}$ one may wonder whether the classes $\mathcal C _{\kappa}$ fit in the setting of flat relative Mittag-Leffler modules also for $\kappa > \aleph_0$.  

It is easy to see that $\mathcal F \mathcal M \subseteq \mathcal C _{\aleph_1} \subseteq  \mathcal F \mathcal P$ for any ring $R$, and that $\mathcal C _{\aleph_1} =  \mathcal F \mathcal M$ when $R$ is right hereditary or von Neumann regular (cf.\ \cite[3.19]{GT} and \cite[3.7(iii)]{BT}). Also, for each $\kappa \geq \aleph_0$, all $< \kappa$-generated modules in the class $\mathcal C _\kappa$ are projective. In particular, for $\kappa = \aleph_1$, the classes $\mathcal C _{\aleph_1}$ and $\mathcal F \mathcal M$ contain the same countably presented modules (namely the projective ones), so if $\mathcal C _{\aleph_1} \neq \mathcal F \mathcal M$, then $\mathcal C _{\aleph_1} \neq \mathcal D _{\mathcal Q}$ for any class of left $R$-modules $\mathcal Q$ by \cite[2.5(i)]{BT}. If $R$ is not right perfect, then the class $\mathcal F \mathcal M$ contains $\aleph_1$-generated non-projective modules (cf.\ \cite[3.19]{GT} and \cite[VII.1.3]{EM}), so for each $\kappa > \aleph_1$, $\mathcal F \mathcal M \nsubseteq \mathcal C _{\kappa}$, whence again $\mathcal C _{\kappa} \neq \mathcal D _{\mathcal Q}$ for any class of left $R$-modules $\mathcal Q$. 
\end{remark}
\medskip

\medskip
\begin{definition} Let $R$ be a commutative ring.
\begin{enumerate}
\item[(1)] Let $\mathfrak{P} $ a property of modules. Then $ \mathfrak{P}(\ModR) $ denotes the class of all modules satisfying the property $ \mathfrak{P} $.
\item[(2)] Let $\mathfrak{R}$ be a class of commutative rings. Let $X$ be a scheme and $(\mathcal R (u)\;|\; u \subseteq X, u \hbox{ open affine })$ be its structure sheaf. Then $X$ is a \emph{locally $\mathfrak{R}$-scheme} provided that $\mathcal R (u) \in \mathfrak{R}$ for each open affine set $u$ of $X$.
\end{enumerate}
\end{definition}

The main properties $\mathfrak{P} $ of modules that we will be interested in here are the flatness, projectivity, and various properties related to Mittag-Leffler conditions that are in general weaker than projectivity, but stronger than flatness.  We will work with general schemes, but in our final application, we will restrict ourselves to \emph{locally noetherian} schemes, that is, the locally $\mathfrak{R}$-schemes where $\mathfrak{R}$ is the class of commutative noetherian rings.    

\medskip  
Recall that given two commutative rings $ R $ and $ S $, a ring homomorphism $ \varphi : R \rightarrow S $ is \emph{flat}, provided that $ S $ is a flat $ R $-module (where the $ R $-module structure on $ S $ is induced by $ \varphi $), that is, the functor $ F = -\otimes_{R} S $ is exact.\\
Moreover, $ \varphi $ is \emph{faithfully flat} provided that $ \varphi $ is flat, and $ N \otimes_{R} S \neq 0$, whenever $ 0 \neq N \in \ModR$. Faithful flatness of $ \varphi $ is equivalent to the following property of the functor $ F $: for each complex $ \mathcal{C} $ of $ R $-modules, $ \mathcal{C} $ is exact in $\ModR$, if and only if $ F(\mathcal{C}) $ is exact in $\ModS$, \cite[Theorem 7.2]{M}. 

A useful characterization of faithfully flat ring homomorphisms of commutative rings goes back to \cite[Chap. I, \S3, Proposition 9]{B} (see also \cite[Lemma 2]{An}):

\begin{lemma}\label{ang1} A flat ring homomorphism $ \varphi : R \rightarrow S $ of commutative rings is faithfully flat, if and only if $ \varphi$ -- viewed as an $R$-homomorphism -- is a pure monomorphism. 
\end{lemma}

Next, we recall the classic notions of ascent and descent, cf.\ \cite[10.82]{St} or \cite{P}. 

\begin{definition}\label{defad}
	Let $ \mathfrak{P} $ be a property of modules, and $\mathfrak{R}$ a class of commutative rings.
	\begin{enumerate}
		\item[(1)] $ \mathfrak{P} $ is said to \emph{ascend} along flat morphisms in $\mathfrak{R}$, provided that for each flat ring homomorphism $ \varphi : R \rightarrow S, $ such that $R, S \in \mathfrak{R}$, and each $ M \in \mathfrak{P}(\ModR) $, also $ M \otimes_{R}S \in \mathfrak{P}(\ModS) $. 
		\item[(2)] $ \mathfrak{P} $ is said to \emph{descend} along faithfully flat morphisms in $\mathfrak{R}$, provided that for each faithfully flat ring homomorphism of commutative rings, $ \varphi : R \rightarrow S $, such that $R \in \mathfrak{R}$ and $S$ is a finite direct product of rings from $\mathfrak{R}$, and for each $ M \in \ModR $, such that $ M \otimes_{R}S \in \mathfrak{P}(\ModS) $, also $ M \in \mathfrak{P}(\ModR) $.
		\item[(3)] $ \mathfrak{P} $ is an \emph{ad-property} in $\mathfrak{R}$ , provided that $ \mathfrak{P} $ ascends along flat morphisms in $\mathfrak{R}$, descends along faithfully flat morphisms in $\mathfrak{R}$, and, moreover, $ \mathfrak{P} $ is \emph{compatible with finite ring direct products} in the following sense: if $ R =\prod_{i<n} R_{i} $ is a finite ring direct product of rings with $R_i \in \mathfrak{R}$ for each $i < n$, and $ (M_{i} \mid i < n) $ satisfy $ M_{i} \in \mathfrak{P}(\ModR_{i}) $ for each $ i<n$, then  $ M =\prod_{i<n} M_{i} \in \mathfrak{P}(\ModR) $.
	\end{enumerate}
\end{definition}

In the case when $\mathfrak{R}$ is the class of all commutative rings, we will omit the attribute {\lq}in $\mathfrak{R}${\rq} and say simply that $ \mathfrak{P} $ ascends, descends, and $ \mathfrak{P} $ is an ad-property.

\medskip
Let $ \mathfrak{P} $ be a property of modules. If $X$ is an affine scheme, i.e., $X = \Spec {R}$ for a commutative ring $R$, then $\Qcoh {X} = \rmod R$, so $ \mathfrak{P} $ is at the same time a property of quasi-coherent sheaves on $X$. For general schemes $X$, one can extend $ \mathfrak{P} $ to a property of quasi-coherent sheaves $\mathcal M$ on $X$ algebraically, by requiring property $ \mathfrak{P} $ to hold for each module of sections of $\mathcal M$:      

\begin{definition}\label{induced}
Let $ \mathfrak{P} $ be a property of $R$-modules, $X$ a scheme, and $(\mathcal R (u)\;|\; u \subseteq X, u \hbox{ open affine })$ be its structure sheaf. A quasi-coherent sheaf $\mathcal M$ on $X$ is a \emph{locally $\mathfrak P$-quasi-coherent sheaf on $X$} in the case when for each open affine set $u$ of $X$, the $\mathcal R (u)$-module of sections $\mathcal M (u)$ satisfies $\mathfrak P$. That is, $\mathcal M (u) \in \mathfrak P (\mathcal R (u))$.
\end{definition}  

If $ \mathfrak{P} $ is the property of being a projective module, then the locally $\mathfrak P$-quasi-coherent sheaves are the (infinite dimensional) vector bundles, see \cite{D}. When $ \mathfrak{P} $ denotes the property of being a flat Mittag-Leffler module (a restricted flat Mittag-Leffler module) then by \cite{EGPT}, the locally $\mathfrak P$-quasi-coherent sheaves are called \emph{Drinfeld vector bundles} (\emph{restricted Drinfeld vector bundles}). Extending this notation to $n \geq 0$, we will call a quasi-coherent sheaf $\mathcal M$ an \emph{$n$-Drinfeld vector bundle} in case it is a locally $\mathfrak P _n$-quasi-coherent sheaf where $\mathfrak P _n$ is the property of being a flat $n$-Mittag-Leffler module. Thus, $0$-Drinfeld vector bundles are just the Drinfeld vector bundles from \cite{EGPT}.  

\medskip
A basic question concerning the various algebraic notions of locally $\mathfrak P$-quasi-coherent sheaves defined above is whether these notions are also geometric, independent on a particular choice of affine coordinates on $X$, that is, whether the notions are Zariski local:

\begin{definition}\label{localprop} Let $\mathfrak R$ be a class of commutative rings, and $\mathfrak C $ be the class of all locally $\mathfrak{R}$-schemes. 

The notion of a locally $\mathfrak P$-quasi-coherent sheaf is \emph{Zariski local on $\mathfrak C$} provided that for each $X \in \mathfrak C$, each open affine covering $X = \bigcup_{v \in V} v$ of $X$, and each quasi-coherent sheaf $\mathcal M$ on $X$, the following implication holds true: if $\mathcal M (v) \in \mathfrak P (\mathcal R (v))$ for all $v \in V$, then $\mathcal M$ is locally $\mathfrak P$-quasi-coherent. 
\end{definition}

ad-properties of modules are important, because they guarantee Zariski locality:

\begin{lemma}\label{adtoZar} Let $\mathfrak R$ be a class of commutative rings. Let $\mathfrak P$ be an ad-property in $\mathfrak R$. Then the notion of a locally $\mathfrak P$-quasi-coherent sheaf is Zariski local on the class of all locally $\mathfrak{R}$-schemes.
\end{lemma}
\begin{proof} This is proved via the Affine Communication Lemma \cite[5.3.2]{V}, see \cite[27.21.2]{St} or \cite[Lemma 3.5]{EGT}.
\end{proof}
   
It is well-known that the properties of being a projective, flat, flat Mittag--Leffler, and restricted flat Mittag--Leffler module, are ad-properties in the class of all commutative rings. Thus the corresponding notions of an (infinite dimensional) vector bundle, flat quasi-coherent sheaf, Drinfeld vector bundle, and restricted Drinfeld vector bundle, are Zariski local on the class of all schemes (see \cite[Seconde partie]{RG}, \cite[\S\S{8-9}]{P}, and \cite{EGT}). Further instances of ad-properties, related to tilting and silting, have recently been introduced in \cite{BHM} and \cite{HST}.  

Our goal here is to investigate the ascent and descent for flat relative Mittag-Leffler modules, i.e., the flat $\mathcal Q$-Mittag-Leffler modules where $\mathcal Q$ is a subclass of $\lmod R$. Then we will apply the results obtained to proving Zariski locality for the corresponding notions of quasi-coherent sheaves. 

\medskip 
For further unexplained terminology, we refer to \cite{EJ} and \cite{GT}.

\section{The algebraic background of ascent and descent for flat relative Mittag-Leffler modules}\label{tech}

First, we recall some connections between the Mittag-Leffler property and stationarity.

\begin{definition}
Let $R$ be an arbitrary ring and $ B $ be a module. 
	\begin{enumerate}
		\item[(1)] Let $( I , \leq )$ be an upper directed poset. A direct system $ (M_{i}, f_{ji} \;|\; i \leq j \in I) $ of modules is said to be $ B $-\emph{stationary} provided that the induced inverse system 
		$$(\Hom R {M_{i}} {B}, \Hom R {f_{ji}} {B} \;|\; i \leq j \in I)$$ 
		satisfies the \emph{Mittag-Leffler condition}, that is, for each $i \in I$ there exists $i \leq j \in I$ such that $\im 
		\Hom R {f_{ki}} {B} = \im \Hom R {f_{ji}} {B}$ for all $j \leq k \in I$.
		\item[(2)] A module $ M $ is said to be $ B $-\emph{stationary} if there exists a $ B $-stationary direct system of finitely presented modules $ (M_{i}, f_{ji} \;|\; i \leq j \in I) $ such that $ M = \varinjlim M_{i} $.
		\item[(3)] Let $ \mathcal B $ be a class of right $ R $-modules. We say that a direct system $ (M_{i}, f_{ji} \;|\; i \leq j \in I) $, or a right $ R $-module $ M $, is $ \mathcal B $-\emph{stationary}, if it is $ B $-stationary for all $ B \in \mathcal B$.
	\end{enumerate}
\end{definition}

Recall that a class of modules is said to be \emph{definable} provided that it is closed under direct limits, direct products and pure submodules. For each class of modules $\mathcal Q$ there is the least definable class of modules containing $\mathcal Q$, called the \emph{definable closure} of $\mathcal Q$ and denoted by $\Def \mathcal Q$. It is obtained by closing $\mathcal Q$ first by direct products, then direct limits, and finally by pure submodules, cf.\ \cite[Lemma 2.9 and Corollary 2.10]{H}. Note that each definable class is also closed under direct sums, pure extensions, and pure-epimorphic images (see e.g.\ \cite[Lemma 6.9]{GT}).

There is a duality between definable classes of left and right $R$-modules: given a definable class $\mathcal Q$ of left (right) $R$-modules, the \emph{dual definable class} $\mathcal Q ^{\vee}$ of $\mathcal Q$ is the least definable class of right (left) $R$-modules containing the character modules $Q^+ = \Hom {\mathbb Z}{M}{\mathbb Q/\mathbb Z}$ of all modules $M \in \mathcal Q$. Then $\mathcal Q = (\mathcal Q ^{\vee})^{\vee}$ for any definable class of left (right) modules $\mathcal Q$, see e.g.\ \cite[\S 2.5]{R2}.

\medskip

FP$_2$ modules are important sources of mutually dual definable classes of left and right modules:  
 
\begin{example}\label{coherent}
Let $\mathcal S$ be a class of FP$_2$ modules. Then $\mathcal S ^\perp$ is a definable class in $\rmod R$ (see \cite[Example 6.10]{GT}), and $\mathcal S ^\intercal$ is a definable class of left $R$-modules. Indeed, $\mathcal S ^\intercal$ is always closed under direct limits and pure submodules, and since $\mathcal S$ consists of FP$_2$ modules, $\mathcal S ^\intercal$ is also closed under products (cf.\ \cite[Theorem 3.2.26]{EJ} and \cite[\S VIII.5]{Br}). Since $M^+ \in \mathcal S ^\intercal$ for each $M \in \mathcal \mathcal S ^\perp$, and $N^+ \in \mathcal S ^\perp$ for each $N \in \mathcal S ^\intercal$ by \cite[Lemma 2.16(b) and (d)]{GT}, the definable classes $\mathcal S ^\perp$ and $\mathcal S ^\intercal$ are mutually dual. 

The classes of left $R$-modules $\mathcal Q$ of the form  $\mathcal Q = \mathcal S ^\intercal$ for a class $\mathcal S$ consisting of FP$_2$ modules will be called \emph{of finite type}.

For example, when $R$ is a right coherent ring and $\mathcal S$ the class of all finitely presented modules, then the class $\mathcal S ^\perp$ of all absolutely pure modules is definable in $\rmod R$, and its dual definable class of all flat left $R$-modules, $\mathcal S ^\intercal$, is of finite type.      
\end{example}

\begin{proposition}\cite[Proposition 1.7 and Theorem 2.11]{H}\label{thm1}
	Let $ R $ be a ring. Let $ \mathcal Q  $ be a definable class of left $R$-modules and $ \mathcal B = \mathcal Q ^{\vee}$ be its dual definable class. Let $ M $ be a right $ R $-module.
	Then the following conditions are equivalent:
	\begin{enumerate}
		\item[(1)] $ M $ is $ \mathcal Q  $-Mittag-Leffler.
		\item[(2)] $ M $ is $ Q $-Mittag-Leffler for all $ Q \in \mathcal Q $.
		\item[(3)] $ M $ is $ Q^{+} $-stationary for all $ Q \in \mathcal Q$.
		\item[(4)] $ M $ is $ \mathcal B $-stationary.
	\end{enumerate}
\end{proposition}

While studying flat $ \mathcal Q $-Mittag-Leffler modules, one can actually restrict to definable classes of modules $\mathcal Q$:

\begin{proposition}\cite[Corollary 2.10]{H}\label{prop1} Let $ \mathcal Q $ be a class of left $ R $-modules. Let $M$ be a $\mathcal Q$-Mittag-Leffler module. Then $M$ is also $\Def \mathcal Q$-Mittag-Leffler. 
\end{proposition}

Now we will turn to the ascent for flat relative Mittag-Leffler modules, so we will again restrict ourselves to commutative rings. 

\begin{lemma}\label{ascent}
Let $ \varphi : R \rightarrow S  $ be a flat homomorphism of commutative rings and $\mathcal Q$ be any class of modules. If $ M $ is a flat $ \mathcal Q $-Mittag-Leffler module, then $ M \otimes_{R}S $ is a flat $(\mathcal Q \otimes_{R} S)$-Mittag-Leffler $ S $-module.
\end{lemma}
\begin{proof}
	Since $M$ is a flat module, the functor $(M \otimes_{R} S) \otimes_S - : \rmod S \to \rmod {\mathbb Z}$ is a composition of two exact functors
	\begin{center}
		$  (M \otimes_{R} S)\otimes_{S}- = (M \otimes_{R} -)(S \otimes_{S}-)$.
	\end{center}
	Thus $ M \otimes_{R}S $ is a flat $S$-module.\\
	Assume that $ M $ is a $ \mathcal Q $-Mittag-Leffler module and let $( Q_i \mid i \in I )$ be a family of elements of $\mathcal Q$. First, note that $ \mathcal Q \otimes_{R} S \subseteq \Def (\mathcal Q)$ as classes of modules. Indeed, since $S$ is a flat module, we can write it as a direct limit of finitely generated free modules, say $ S= \varinjlim_{\alpha} R^{n_{\alpha}} $. Therefore, $ \mathcal Q \otimes_{R} \varinjlim_{\alpha} R^{n_{\alpha}} \cong \varinjlim_{\alpha} \mathcal Q^{n_{\alpha}} \in \Def (\mathcal Q )$. 	By our assumption on $M$ and by Proposition \ref{prop1}, we infer that the canonical map $\psi_M : M \otimes_{R} {\prod_{i \in I} (Q_i \otimes_{R} S)} \to \prod_{i \in I} (M\otimes_{R} Q_i \otimes_{R} S)$ is monic.

We have the following commutative diagram whose horizontal maps are isomorphisms:
$$\begin{CD}
   {(M\otimes_{R} S) \otimes_{S} {\prod_{i \in I} (Q_i \otimes_{R} S)}} @>{\cong}>>  {M \otimes_{R} {\prod_{i \in I} (Q_i \otimes_{R} S)}}\\ 
     @V{\psi_{M \otimes_R S}}VV  @V{\psi_M}VV \\
      {\prod_{i \in I} (M \otimes_{R} S) \otimes_{S}(Q_i \otimes_{R} S)} @>{\cong}>> {\prod_{i \in I}  M\otimes_{R} (Q_i \otimes_{R} S)}. \\
\end{CD}$$
Here, the left vertical map $\psi_{M \otimes_R S}$ is the canonical morphism $\psi_{M \otimes_R S} : (M\otimes_{R} S) \otimes_{S} {\prod_{i \in I} (Q_i \otimes_{R} S)} \to \prod_{i \in I} (M\otimes_{R} S) \otimes_{S} (Q_i \otimes_{R} S)$. Thus $\psi_{M \otimes_R S}$ is monic. This proves that  $ M \otimes_{R} S $ is a $ (\mathcal Q \otimes_{R} S)$-Mittag-Leffler $ S $-module.
\end{proof}

The descent of flatness is well-known, we include a proof here for the sake of completeness.

\begin{lemma}\label{known}
	Let $ \varphi : R \rightarrow S  $ be a faithfully flat homomorphism of commutative rings, and let $ M $ be a module such that the $ S $-module $ M \otimes_{R}S $ is flat. Then $ M $ is a flat module.
\end{lemma}
\begin{proof}
First, since $S$ is a flat module, also $ M \otimes_{R} S $, viewed as an $R$-module, is flat. Indeed, the functor  $(M \otimes_{R} S) \otimes_R -$ is a composition of two exact functors as follows: $ M \otimes_{R}(S \otimes_{S}S) \otimes_R - = ((M \otimes_{R} S)\otimes_{S} - )(S \otimes_R -) $. So for each short exact sequence $ \mathcal{C} $ of modules, $ \mathcal{C} \otimes_{R}(M \otimes_{R} S) $ is a short exact sequence of $ S $-modules. Hence, by faithful flatness of $ \varphi $, $ \mathcal{C} \otimes_{R} M $ is exact in $ \ModR $, whence $ M $ is a flat module. 
\end{proof}

Recently, a short proof of the descent of the (absolute) flat Mittag-Leffler property along all pure (and hence all faithfully flat) ring homomorphisms was presented in \cite[Lemma 5]{An}.  We include this short proof here as it works also in our relative setting. (We refer to \cite{H2} for a broader context and further applications.)

\begin{lemma}\label{ang2}
Let $ \varphi : R \rightarrow S  $ be a pure monomorphism of commutative rings. Let $\mathcal Q$ be a class of modules. Let $ M $ be a flat module such that  $M \otimes_{R}S $ is a $(\mathcal Q \otimes_R S)$-Mittag-Leffler $S$-module. Then $ M $ is a $\mathcal Q$-Mittag-Leffler module.
\end{lemma}
\begin{proof}  
Let $( Q_i \mid i \in I )$ be a family consisting of modules from $\mathcal Q$. Since $\varphi$ is pure, the canonical morphism $g_i : Q_i \cong Q_i \otimes_R R \to Q_i \otimes_R S$ is monic for each $i \in I$, and so is $g = \prod_{i \in I} g_i : \prod_{i \in I} Q_i \to \prod_{i \in I} (Q_i \otimes_R S)$. 

Let $M$ be a flat module such that  $M \otimes_{R}S $ is a $(\mathcal Q \otimes_R S)$-Mittag-Leffler $S$-module. Since $M$ is flat, the morphism $M \otimes_R g : M \otimes_R \prod_{i \in I} Q_i \to M \otimes_R \prod_{i \in I} (Q_i \otimes_R S)$ is monic. Moreover, we have the canonical isomorphism $\psi : M \otimes_R \prod_{i \in I} (Q_i \otimes_R S) \cong M \otimes_R (S \otimes_S \prod_{i \in I} (Q_i \otimes_R S)) \cong (M \otimes_R S) \otimes_S \prod_{i \in I} (Q_i \otimes_R S)$. Since $M \otimes_{R} S $ is a $(\mathcal Q \otimes_R S)$-Mittag-Leffler $S$-module, the canonical morphism $h : (M \otimes_R S) \otimes_S \prod_{i \in I} (Q_i \otimes_R S) \to \prod_{i \in I} (M \otimes_R S) \otimes_S (Q_i \otimes_R S)$ is monic. Thus the composite morphism $k = h \psi (M \otimes_R g)$ is monic. 

Notice that $k (m \otimes_R (q_i)_{i \in I}) = ((m \otimes_R 1) \otimes_S ( q_i \otimes_R 1))_{i \in I}$, so $k$ can also be expressed as the composition of another triple of canonical morphisms: 
$k = \psi^\prime g^\prime h^\prime$, where $h^\prime : M \otimes_R \prod_{i \in I} Q_i \to \prod_{i \in I} (M \otimes_R Q_i)$, $g^\prime$ is the monomorphism $ \prod_{i \in I} (M \otimes_R Q_i) \to \prod_{i \in I} (M \otimes_R Q_i \otimes_R S)$, and $\psi ^\prime$ the isomorphism $\prod_{i \in I} (M \otimes_R Q_i \otimes_R S) \to \prod_{i \in I} (M \otimes_R S) \otimes_S (Q_i \otimes_R S)$. Since $k$ is monic, so is $h^\prime$. The latter says that $ M $ is a $\mathcal Q$-Mittag-Leffler module.
\end{proof}

Now, we can easily prove the descent for flat relative Mittag-Leffler modules:        

\begin{theorem}\label{descent}
Let $ \varphi : R \rightarrow S  $ be a faithfully flat homomorphism of commutative rings. Let $\mathcal Q$ be a class of modules. Let $M$ be a module such that $M \otimes_R S$ is a flat $(\mathcal Q \otimes_R S)$-Mittag-Leffler $S$-module. Then $M$ is a flat $\mathcal Q$-Mittag-Leffler module. 
\end{theorem} 
\begin{proof} By Lemma \ref{known}, we can assume that $M$ is a flat module. By Lemma \ref{ang1}, $ \varphi$ is a pure monomorphism, so Lemma \ref{ang2} applies and shows that $ M $ is a $\mathcal Q$-Mittag-Leffler module.  
\end{proof}

\medskip
It is worth noting that for countably presented flat modules, Mittag-Leffler conditions relative to definable classes of modules can be expressed in terms of vanishing of the Ext functor, following \cite[\S1]{H}. 

\begin{lemma}\label{herb} Let $R$ be any ring. Let $M$ be a countably presented flat module, $\mathcal Q$ be a definable class of left $R$-modules, and $ \mathcal B = \mathcal Q ^{\vee}$. Then $M$ is $\mathcal Q$-Mittag-Leffler, if and only if $M \in {}^\perp \mathcal B$. 
\end{lemma}
\begin{proof} If $M$ is $\mathcal Q$-Mittag-Leffler, then $M$ is $\mathcal B$-stationary by Proposition \ref{thm1}. Since $M$ is a countable direct limit of finitely presented free modules and $\mathcal B$ is closed under countable direct sums, we infer from \cite[Corollary 2.23]{GT} and \cite[Lemma 1.11(3)]{H} that $\Ext 1RMB = 0$ for each $B \in \mathcal B$. The converse implication follows by \cite[Lemma 1.11(1)]{H} and Proposition \ref{thm1}. 
\end{proof}

\begin{remark} Lemma \ref{herb} does not extend to uncountably presented modules in general. Just consider any non-right perfect ring $R$ and let $\mathcal Q = \lmod R$. Then $ \mathcal B = \mathcal Q ^{\vee} = \rmod R$, so ${}^\perp \mathcal B = \mathcal P_0 \subsetneq \mathcal F \mathcal M$ (though, as correctly claimed by Lemma \ref{herb}, the countably presented modules in $\mathcal P_0$ and $\mathcal F \mathcal M$ are the same).
\end{remark}   

Theorem \ref{descent}, Proposition \ref{prop1}, and Lemmas \ref{ascent} and \ref{herb} yield the following 

\begin{corollary}\label{corcount}
Let $ \varphi : R \rightarrow S  $ be a faithfully flat homomorphism of commutative rings. Let $\mathcal Q$ be a definable class of modules and $ \mathcal B = \mathcal Q ^{\vee}$. Let $\mathcal Q ^\prime$ denote the least definable class of $S$-modules containing $\mathcal Q \otimes_R S$, and $\mathcal B ^\prime$ its dual definable class.    

Let $M$ be a countably presented flat module. Then $M \in {}^\perp \mathcal B$, if and only if $M \otimes_R S \in {}^\perp \mathcal B ^\prime$.   
\end{corollary}

\medskip

In the particular setting of definable classes arising from kernels of Tor functors (such as the definable classes of finite type from Example \ref{coherent}), we have the following relation between definable closures:

\begin{lemma}\label{least} Let $ \varphi : R \rightarrow S  $ is a flat homomorphism of commutative rings and $\mathcal C$ be a class of $R$-modules. Then $\Def {(\mathcal C \otimes_R S)^\intercal} = \Def {(\mathcal C ^\intercal \otimes_R S )}$.

In particular, if $\mathcal C$ consists of FP$_2$ modules, then $\Def {(\mathcal C ^\intercal \otimes_R S )} = (\mathcal C \otimes_R S)^\intercal$.
\end{lemma} 
\begin{proof} First, $(\mathcal C ^\intercal) \otimes_R S \subseteq  (\mathcal C \otimes_R S)^\intercal$ by Lemma \ref{Ext-Tor}(1), whence 
$\Def {(\mathcal C ^\intercal \otimes_R S )} \subseteq \Def {(\mathcal C \otimes_R S)^\intercal}$. 

For the opposite inclusion, note that by Lemma \ref{Ext-Tor}(2), $(\mathcal C \otimes_R S)^\intercal$ is the class of all $S$-modules $N$ satisfying the following condition: $N$, viewed as an $R$-module, is an element of $\mathcal C ^\intercal$. Then again $N \otimes_R S \in (\mathcal C \otimes_R S)^\intercal$ by Lemma \ref{Ext-Tor}(1). Since the canonical homomorphism $f : n \mapsto n \otimes 1$ from $N$ to $N \otimes_R S$ is an $S$-homomorphism, and the $S$-homomorphism $g : N \otimes_R S \to N$ defined by $g : n \otimes s \mapsto n.s$ satisfies $g f = 1_N$, we infer that $N$ is isomorphic to a direct summand in $N \otimes_R S$ as an $S$-module. Thus $(\mathcal C \otimes_R S)^\intercal$ consists of $S$-modules isomorphic to direct summands of the modules from $\mathcal C ^\intercal \otimes_R S$, whence $(\mathcal C \otimes_R S)^\intercal \subseteq \Def {(\mathcal C ^\intercal \otimes_R S )}$, proving the opposite inclusion.   

If $\mathcal C$ consists of FP$_2$ modules, then also $\mathcal C \otimes_R S$ consists of FP$_2$ $S$-modules, whence $(\mathcal C \otimes_R S)^\intercal$ is a definable class by Example \ref{coherent}. 
\end{proof}

\section{Zariski locality of quasi-coherent sheaves associated with flat relative Mittag-Leffler modules}\label{appl}

In this section, we will apply the results of Section \ref{tech} to prove Zariski locality of flat relative $\mathcal Q$-Mittag-Leffler modules in various particular settings. 

We start with a direct general application to quasi-coherent sheaves associated with f-projective modules. Recall that a module $M$ is \emph{f-projective} if $M$ is flat and $\{ R \}$-Mittag-Leffler, or equivalently, $M$ is a flat $\mathcal Q$-Mittag-Leffler module where $\mathcal Q$ is the class of all flat left $R$-modules, \cite{G} (see also Proposition \ref{prop1} and \cite[\S3]{BT}). In accordance with our Definition \ref{induced}, we call a quasi-coherent sheaf $\mathcal M$ on a scheme $X$ \emph{locally f-projective} in case for each open affine set $u$ in $X$, the $\mathcal R (u)$-module of sections $\mathcal M (u)$ is an f-projective $\mathcal R (u)$-module.

\begin{theorem}\label{f-locality} The notion of a locally f-projective quasi-coherent sheaf is Zariski local on the class of all schemes. 
\end{theorem}
\begin{proof} By Lemma \ref{adtoZar}, it suffices to prove that the property of being an f-projective module is an ad-property in the class of all commutative rings. However, its ascent and descent follows for $\mathcal Q = \{ R \}$ immediately by Lemma \ref{ascent} and Theorem \ref{descent}, respectively. The compatibility with finite ring direct products is obvious (cf.\ Definition \ref{defad}(3)).    
\end{proof}

\medskip

\emph{For the rest of this section, $R$ will denote a commutative ring, $\mathcal C _{R}$ a class of modules, and $\mathcal Q _{R}$ the definable class $\mathcal Q _{R} = \Def {\mathcal C _{R}^\intercal}$.} In particular, $\mathcal Q _{R} = \mathcal C _{R}^\intercal$ in case $\mathcal C _{R}$ consists of FP$_2$ modules.

The relevant property $\mathfrak{P}$ of modules is defined as follows: \emph{if $M$ is a module, then $M \in \mathfrak{P}(\rmod R)$, iff $M$ is a flat $\mathcal Q _{R}$-Mittag-Leffler module.}    

\medskip

In order to prove locality of the induced notions of quasi-coherent sheaves in this setting, we will need compatibility of the properties $\mathfrak{P}$ for commutative rings $R$ and $S$ connected by flat, and faithfully flat, morphisms. More precisely, we will require the following compatibility conditions (C1), (C2) and (C3):   

\begin{definition} Let $\mathfrak{R}$ be a class of commutative rings.
\begin{enumerate}
		\item[(C1)] For each flat ring homomorphism $ \varphi : R \rightarrow S  $ with $R, S \in \mathfrak R$, $\mathcal C _{R} \otimes_R S \subseteq \mathcal C _{S}$.
		\item[(C2)] For each faithfully flat ring homomorphism $ \varphi : R \rightarrow S  $ where $R \in \mathfrak R$ and $S$ is a finite direct product of rings in $\mathfrak R$, $\Def {\mathcal C _{S}^\intercal} = \Def {(\mathcal C _{R} \otimes_R S)^\intercal}$. 
		\item[(C3)] If $S = \prod_{i < n} R_i$ where $R_i \in \mathfrak R$ for each $i < n$, then $\mathcal C _{S} = \prod_{i < n} \mathcal C _{R_i}$.
\end{enumerate}
Notice that (C1) implies the inclusion $\mathcal C _{S}^\intercal \subseteq (\mathcal C _{R} \otimes_R S)^\intercal$, and hence $\Def {\mathcal C _{S} ^\intercal} \subseteq \Def {(\mathcal C _{R} \otimes_R S)^\intercal}$. 
\end{definition} 

\begin{lemma}\label{ascents} Let $\mathfrak R$ be a class of commutative rings such that condition (C1) holds. Then the property $ \mathfrak{P} $ ascends along flat morphisms in $\mathfrak{R}$.
\end{lemma}
\begin{proof} Let $ \varphi : R \rightarrow S  $ be a flat ring homomorphism with $R, S \in \mathfrak R$ and $M$ be a flat $\mathcal Q _{R}$-Mittag-Leffler module. By Lemma \ref{ascent}, $M \otimes_R S$ is a flat $(\mathcal Q _{R} \otimes_R S)$-Mittag-Leffler $S$-module, and hence a flat $\Def (\mathcal Q _{R} \otimes_R S )$-Mittag-Leffler $S$-module by Proposition \ref{prop1}. Condition (C1) and Lemma \ref{least} give 
$$\mathcal Q _{S} = \Def {\mathcal C _{S} ^\intercal} \subseteq \Def {(\mathcal C _{R} \otimes_R S)^\intercal} = \Def (\mathcal Q _{R} \otimes_R S ).$$ Thus, $M \otimes_R S$ is a flat $\mathcal Q _{S}$-Mittag-Leffler $S$-module. 
\end{proof}

\begin{lemma}\label{descents} Let $\mathfrak R$ be a class of commutative rings such that condition (C2) holds. Then the property $\mathfrak{P}$ descends along faithfully flat morphisms in $\mathfrak{R}$.
\end{lemma}
\begin{proof} Let $ \varphi : R \rightarrow S  $ be a faithfully flat ring homomorphism, where $R \in \mathfrak R$ and $S$ is a finite direct product of rings in $\mathfrak R$. Let $M$ be a module such that $M \otimes_R S$ is a flat $\mathcal Q _{S}$-Mittag-Leffler $S$-module. Condition (C2) and Lemma \ref{least} yield 
$$\mathcal Q _{S} = \Def {\mathcal C _{S}^\intercal} = \Def {(\mathcal C _{R} \otimes_R S)^\intercal} = \Def (\mathcal Q _{R} \otimes_R S ),$$ so $M \otimes_R S$ is a flat $(\mathcal Q _{R} \otimes_R S)$-Mittag-Leffler $S$-module. By Theorem \ref{descent}, $M$ is a flat $\mathcal Q _{R}$-Mittag-Leffler module. 
\end{proof}

Thus, we obtain 

\begin{theorem}\label{locality} Let $\mathfrak R$ be a class of commutative rings such that conditions (C1), (C2) and (C3) hold.
Then $\mathfrak{P}$ is an ad-property in $\mathfrak R$, whence the notion of a locally $\mathfrak{P}$-quasi-coherent sheaf is Zariski local on the class of all locally $\mathfrak R$-schemes.
\end{theorem}
\begin{proof} By condition (C3), $\mathfrak{P}$ is compatible with finite ring direct products, so the ad-property of $\mathfrak{P}$ follows by Lemmas \ref{ascents} and \ref{descents}. The final claim follows by Lemma \ref{adtoZar}.  
\end{proof}

We finish this section by noting several applications of Theorem \ref{locality}:

\medskip
\subsection{Applications}
\medskip

{\bf 1.} Let $\mathcal R$ be the class of all commutative rings and $\mathcal C _R = \{ 0 \}$, so $\mathcal Q _R = \lmod R$. In this case, Theorem \ref{locality} yields the Zariski locality of the notion of a Drinfeld vector bundle (= locally flat Mittag-Leffler quasi-coherent sheaf) proved in \cite{EGT}. 

\medskip
{\bf 2.} Let $\mathfrak R$ be the class of all commutative rings and $\mathcal C _R$ the class of all finitely presented modules. Then $\mathcal Q _{R} = \Def {\mathcal C _{R}^\intercal} = 
\Def {\mathcal F _0}$. By Proposition \ref{prop1}, a module $M$ has property $\mathfrak P$, iff $M$ is f-projective. Conditions (C1) and (C3) clearly hold true.  

Condition (C2) holds even in the stronger form of $\mathcal C _{S} ^\intercal = (\mathcal C _{R} \otimes_R S)^\intercal$ whenever $ \varphi : R \rightarrow S  $ is a faithfully flat homomorphism of commutative rings. Indeed, $\mathcal C _{S} ^\intercal$ is the class of all flat $S$-modules. Let $M \in (\mathcal C _{R} \otimes_R S)^\intercal$. By Lemma \ref{Ext-Tor}(2), $\Tor 1R{\mathcal C _{R}}M = 0$, whence $M$ is a flat $R$-module. Then $M \otimes_R S$ is a flat $S$-module, by (the proof of) Lemma \ref{ascent}. However, the $S$-module $M$ is isomorphic to a direct summand in $M \otimes_R S$ (cf.\ the proof of Lemma \ref{least}), whence $M$ is a flat $S$-module. This proves the inclusion $(\mathcal C _{R} \otimes_R S)^\intercal \subseteq \mathcal C _{S} ^\intercal$; the other inclusion is a consequence of condition (C1). 

Thus, Theorem \ref{f-locality} is just a particular instance of Theorem \ref{locality} for $\mathcal C _R$ = the class of all finitely presented modules. 

\medskip
{\bf 3.} A more involved application of Theorem \ref{locality} concerns the case when $\mathcal{C}_R = \mathcal F _n$ for some $n \geq 1$. In this case, we will verify conditions (C1) -- (C3) for $\mathfrak R$ = the class of all noetherian rings.

Condition (C1) holds since $\mathcal C _{R} \otimes_R S \subseteq \mathcal C _{S}$ when $S$ is a flat module, and (C3) is obvious. As in Application 2, it only suffices to prove the inclusion $(\mathcal C _{R} \otimes_R S)^\intercal \subseteq \mathcal C _{S}^\intercal$ for each faithfully flat homomorphism of commutative noetherian rings $\varphi : R \rightarrow S$. 

Recall that for an $S$-module $M$, $M^+ = \Hom{\mathbb Z}{M}{\mathbb Q/\mathbb Z}$ denotes the $S$-module of characters of $M$, and for a class of $S$-modules $\mathcal E$, $\mathcal E ^+ = \{ M^+ \mid M \in \mathcal E \}$. We claim that $\mathcal F _n ^+ = \mathcal I _n \cap (\rmod S )^+$. Since character modules of flat modules are injective, the $\subseteq$ inclusion holds. Conversely, let $N = M^+ \in \mathcal I _n$. Since $S$ is noetherian, character modules of injective modules are flat (e.g., by \cite[Lemma 2.16(d)]{GT}), so $N^+ = M^{++} \in \mathcal F _n$. As the class $\mathcal F _n$ is closed under pure submodules and the embedding $M \hookrightarrow M^{++}$ is pure, $M \in \mathcal F _n$ and the claim is proved. Using \cite[Lemma 2.16(b)]{GT} and the fact that the pure embedding $M \hookrightarrow M^{++}$ splits for any pure-injective module $M$, we get $(\mathcal F _n)^\intercal = {}^\perp (\mathcal F _n ^+) = {}^\perp (\mathcal I _n \cap \mathcal P \mathcal I)$. 

Let $\mathcal S _{n,R}$ denote the class of all finitely generated modules that appear as $n$th syzygies in some projective resolution, $\mathcal P$, of a finitely generated module such that $\mathcal P$ consists of finitely generated modules. Then $R \in \mathcal S _{n,R}$, and since $R$ is noetherian, $\mathcal S _{n,R} ^\perp = \mathcal I _n$ by the Baer Test of Injectivity and by dimension shifting. Let $\mathcal D _{n,R}$ denote the class of all modules $M$ that are isomorphic to direct summands of finite extensions of the modules from $\mathcal S _{n,R}$. Then the class $\mathcal D _{n,R}$ is closed under extensions, direct summands, and contains $R$. Moreover, $\mathcal D _{n,R} ^\perp = \mathcal I _n$. By Lemma \ref{closure}, $\varinjlim \mathcal D _{n,R} = {}^\perp (\mathcal I_n \cap \mathcal P \mathcal I )$.  

Finally, let $M \in (\mathcal C _{R} \otimes_R S)^\intercal$. Then Lemma \ref{Ext-Tor}(2) gives $\Tor 1R{\mathcal C _{R}}M = 0$. By the above, $M$, viewed as an $R$-module, is an element of $\varinjlim \mathcal D _{n,R}$. Since $S$ is a flat module, $\mathcal S _{n,R} \otimes _R S \subseteq \mathcal S _{n,S}$, whence also $\mathcal D _{n,R} \otimes_R S \subseteq \mathcal D _{n,S}$. Moreover, the tensor product commutes with direct limits, so $M \otimes_R S \in \varinjlim \mathcal D _{n,S} = \mathcal C _{S} ^\intercal$. As $M$ is isomorphic to a direct summand in $M \otimes_R S$ as an $S$-module, also $M \in \mathcal C _{S} ^\intercal$, and the inclusion $(\mathcal C _{R} \otimes_R S)^\intercal \subseteq \mathcal C _{S}^\intercal$ is proved.   

Recall that if $n \geq 0$ and $\mathcal Q _n = (\mathcal F _n)^\intercal$, then the flat $\mathcal Q _n$-Mittag-Leffler modules are called flat $n$-Mittag-Leffler, and the corresponding quasi-coherent sheaves are the $n$-Drinfeld vector bundles. Thus, we have the following consequence of Theorem \ref{locality} for $\mathcal{C}_R = \mathcal F _n$:

\begin{theorem}\label{fi} For each $n \geq 1$, the notion of an $n$-Drinfeld vector bundle is Zariski local on the class of all locally noetherian schemes.      
\end{theorem}

\begin{remark} If $R$ is a non-right perfect ring (e.g., a commutative noetherian ring of Krull dimension $\geq 1$), then there is a gap between the classes $\mathcal F \mathcal M$ of all flat Mittag-Leffler modules and $\mathcal F$ of all flat modules. In fact, for each class $\mathcal Q$ of left $R$-modules we have 
$\mathcal{F}\mathcal{M} \subseteq \mathcal D _{\mathcal Q} \subseteq \mathcal F$.
Since $\mathcal D _{\mathcal Q} = \mathcal D_ {\Def (\mathcal Q)}$ by Proposition \ref{prop1} and there is only a set of definable classes of modules, there is also only a set of such intermediate classes $\mathcal D _{\mathcal Q}$ between $\mathcal{F}\mathcal{M}$ and $\mathcal F$ (see also \cite[Theorem 3.5(i)]{BT}). 

Of course, the variety of classes of modules between $\mathcal F \mathcal M$ and $\mathcal F$ translates directly into the same variety of classes of locally $\mathfrak P$-quasi-coherent sheaves in the class of all flat quasi-coherent sheaves on the affine scheme $X = \Spec R$, where $R$ is any commutative non-perfect ring (since in this case, $\Qcoh {X}$ is equivalent to $\rmod R$). Moreover, all these classes contain a flat generator, as they contain all vector bundles on $X$.

However, the picture for non-affine schemes may be different, depending on further properties of the schemes. For example, by \cite{SS}, if $X$ is a quasi-compact and quasi-separated scheme, then $\Qcoh {X}$ contains a flat generator, if and only if $X$ is semiseparated (i.e., the intersection of any two open affine sets is affine).      
\end{remark}

\end{document}